\documentclass[11pt,reqno]{amsart}
\usepackage[notref,notcite]{}
\usepackage{graphicx}
\usepackage{amssymb}
\usepackage{amsmath}

\begin{document}

\def\?{
?\vadjust{\vbox to 0pt{\vss\hbox{\kern\hsize\kern1em\large\bf ?!}}}}

\title[]{On a new incomplete Ricci-flat metric.}

\author{E.~G.~Malkovich}
\address{Sobolev Institute of Mathematics and Novosibirsk State University,
         Russia}
\email{malkovich@math.nsc.ru}

\thanks{The author is supported by
a~Grant of the Russian Federation for the State Support of Researches
(Contract No.~14.B25.31.0029).}


\begin{abstract}

We define a system of ODE that gives Einstein 4-dimensional metrics. We found new Ricci-flat incomplete metric of cohomogeneity 1
in explicit formulas and study its characteristics at infinity and in the singular point.



\end{abstract}

\maketitle

\sloppy

\section{Introduction}


One of the main problems in differential geometry is a problem of constructing Ricci-flat and Einstein metrics. Many examples of this metrics were found in the works of the physicists aiming to find solutions for the General Relativity with Lorentzian signature. Still in the Riemannian setting the Eistein equation is very non-trivial and interesting problem to solve.

One can use the Eguchi-Hanson metric to glue K3-surface from the factorized torus by cutting out 16 singular balls and the gluing 16 copies of the Eguchi-Hanson spaces \cite{Donaldson}. The Taub-NUT metric is a Lorentzian Ricci-flat metric in the empty space with topology $\mathbb{R}\times S^3$. The Fubini-Study metric on $\mathbb{C}P^n$ is also well-known and important in differential geometry, for example, as a border case in Sphere Theorem and in many other cases. All this metrics can be considered as a solutions of the system
$$
\left\{\begin{array}{ll}
A_1'=k_1(\frac{A_1}{A_2})^2 +k_2\frac{A_1}{A_2} +k_3, \\
A_2'=l_1(\frac{A_1}{A_2})^2 +l_2\frac{A_1}{A_2} +l_3,
\end{array}\right.
$$
for some definite set of parameters $(k_1,\ldots,l_3)$. We found all parameters for which metric
$$
g=dt^2+A_1^2(t)(e^1)^2+A_2^2(t)((e^2)^2+(e^3)^2)
$$
is Ricci-flat or Einstein and proved the following theorem.

{\bf Theorem 1.}
{\it If metric g is Ricci-flat and functions $A_1,~A_2$ of the metric satisfy the considered system and do not coincide then one of the following cases hold

1. $k_1=1,~k_2=k_3=l_1=0,~ l_2=-1,~ l_3=2$ and metric $g$ is isometric to Taub-NUT metric;

2. $k_1=1,~k_2=0,~k_3=-2,~ l_1=0,~ l_2=-1,~ l_3=0$ and metric $g$ is isometric to Eguchi-Hanson metric;

3. $k_1=1,~k_2=k_3=l_1=0,~ l_2=-1,~ l_3=-2$.}

The case 3 corresponds to a Ricci-flat incomplete metric
$$
\frac{\rho^4}{(c^2-\rho^2)^4}d\rho^2 +\frac{1}{\rho^2}(e^1)^2 +\frac{\rho^2}{(c^2-\rho^2)^2}((e^2)^2+(e^3)^2).
$$
We study asymptotic of this metric at infinity and at singular time. As $t\rightarrow \infty$ the functions $\frac{A_1(t)}{t}\rightarrow 0$ and $A_2(t)\sim 2t$. 
At  singular time $t_0$ the second function $A_2(t)$ looks like  $\gamma (t-t_0)^{\frac{1}{3}}  + \frac{9}{5}(t-t_0)$  for some constant $\gamma$ and the first one $A_1(t)$ blows like  $\frac{\gamma^2}{3}(t-t_0)^{-\frac{1}{3}}$.

We show analogous theorem for Einstein metrics and prove that only Fubini-Study (or hyperbolic Fubini-Study) metric appears in considering setting.

\section{System of ODE for solving Einstein equation.}


We consider standard 3-dimensional sphere $S^3$ with basis of 1-forms $\{e^1,~e^2,~e^3\}$ such that
$$
de^i=2e^{i+1}\wedge e^{i+2}, \quad i=1,2,3 ~mod~ 3.
$$
Then we consider 4-dimensional metric
$$
g=dt^2+A_1^2(t)(e^1)^2+A_2^2(t)((e^2)^2+(e^3)^2)=dt^2+\bar{g}(t)
$$
with orthonormal frame $\varepsilon_0=dt,~ \varepsilon_1=A_1(t)e^1,~ \varepsilon_2=A_2(t)e^2,~ \varepsilon_3=A_2(t)e^3$. For metric $g$ we can calculate connection 1-form $\omega_i^j$ and curvature 2-form $\Omega_i^j$:
$$
-\omega_i^j=\left(  \begin{array}{cccc}
                              0 & \frac{A_1'}{A_1}\varepsilon^1 & \frac{A_2'}{A_2}\varepsilon^2 & \frac{A_2'}{A_2}\varepsilon^3 \\
-\frac{A_1'}{A_1}\varepsilon^1 & 0 & \frac{-A_1}{A_2^2}\varepsilon^3 & \frac{A_1}{A_2^2}\varepsilon^2 \\
-\frac{A_2'}{A_2}\varepsilon^2 & \frac{A_1}{A_2^2}\varepsilon^3 & 0 & \frac{A_1^2-2A_2^2}{A_1A_2^2}\varepsilon^1 \\
-\frac{A_2'}{A_2}\varepsilon^3 & \frac{-A_1}{A_2^2}\varepsilon^2 & -\frac{A_1^2-2A_2^2}{A_1A_2^2}\varepsilon^1 & 0 \\
                            \end{array}
\right)
$$
$$
\Omega^0_1=\varepsilon^0\wedge\varepsilon^1[-\frac{A_1''}{A_1}]+ \varepsilon^2\wedge\varepsilon^3[-\frac{2A_1'}{A_2^2} +\frac{2A_1A_2'}{A_2^3}],
$$
$$
\Omega^0_2=\varepsilon^0\wedge\varepsilon^2[-\frac{A_2''}{A_2}]+ \varepsilon^3\wedge\varepsilon^1[\frac{A_1'}{A_2^2} -\frac{A_1A_2'}{A_2^3}],
$$
$$
\Omega^1_2=\varepsilon^0\wedge\varepsilon^3[\frac{A_1'}{A_2^2}-\frac{A_1A_2'}{A_2^3}]+ \varepsilon^1\wedge\varepsilon^2[\frac{A_1^2}{A_2^4}-\frac{A_1'A_2'}{A_1A_2}],
$$
$$
\Omega^2_3=\varepsilon^0\wedge\varepsilon^1[\frac{2A_1A_2'}{A_2^3}-\frac{2A_1'}{A_2^2}]+ \varepsilon^2\wedge\varepsilon^3[\frac{4}{A_2^2}-\frac{3A_1^2}{A_2^4}-\frac{A_2'^2}{A_2^2}].
$$
We remind that connection form $\omega$ is just the Christoffel symbols rewritten in terms of 1-forms instead of vector fields:
$$
d\varepsilon^i=-\omega^i_j\wedge\varepsilon^j.
$$

Here we use the following action of the differential generated by the standard action on $S^3$:
$$
d\varepsilon^0=ddt=0,\quad d\varepsilon^1=d(A_1e^1)=A_1'dt\wedge\varepsilon^1+2A_1e^2\wedge e^3= \frac{A_1'}{A_1}\varepsilon^0\wedge\varepsilon^1 +  2\frac{A_1}{A_2^2}\varepsilon^2\wedge \varepsilon^3,
$$
$$
d\varepsilon^2=A_2'dt\wedge\varepsilon^2+2A_2e^3\wedge e^1= \frac{A_2'}{A_2}\varepsilon^0\wedge\varepsilon^2 +  2\frac{1}{A_1}\varepsilon^3\wedge \varepsilon^1
$$
$$
d\varepsilon^3=\frac{A_2'}{A_2}\varepsilon^0\wedge\varepsilon^3 +  2\frac{1}{A_1}\varepsilon^1\wedge \varepsilon^2
$$

The curvature form
$$\Omega^i_j =d\omega^i_j +\omega^i_k \wedge \omega^k_j$$
is analog of Riemann tensor $\Omega ^i_j=\frac{1}{2}R^i_{jkl}\varepsilon^k\wedge \varepsilon^l$.

Then one can easily calculate the components of Ricci tensor for the metric $g$ in the frame $\{\varepsilon^0,\ldots \varepsilon^3\}$. They are the following
$$
Ric_{00}=-2\frac{A_1''}{A_1}-4\frac{A_2''}{A_2},\quad Ric_{11}=-2\frac{A_1''}{A_1}- 4\frac{A_1'A_2'}{A_1A_2} +4\frac{A_1^2}{A_2^4},
$$
$$
Ric_{22}=Ric_{33}=-2\frac{A_2''}{A_2} -2\frac{A_1'A_2'}{A_1A_2} -4\frac{A_1^2}{A_2^4}-2\frac{(A_2')^2}{A_2^2} +\frac{8}{A_2^2}.
$$

We notice that for well-known metrics, such as Fubini-Study, Eguchi-Hanson and Taub-NUT, the derivatives of functions $A_1$ and $A_2$ are the polynomials of the variable $\frac{A_1}{A_2}$:
$$
\left\{\begin{array}{ll}
A_1'=k_1(\frac{A_1}{A_2})^2 +k_2\frac{A_1}{A_2} +k_3, \\
A_2'=l_1(\frac{A_1}{A_2})^2 +l_2\frac{A_1}{A_2} +l_3,
\end{array}\right.
\eqno{(1)}
$$
For example, to obtain Fubini-Study metric one should set $k_2=l_1=l_3=0$ and $k_1=2,~k_3=-1,~l_2=1$. We have investigated all possible cases in order to find Ricci-flat and Einstein metrics and prooved

\medskip
{\bf Theorem 1.}
{\it If metric $g$ is Ricci-flat and functions $A_1,~A_2$ of the metric satisfies $(1)$ and do not coincide then one of the following cases hold

1. $k_1=1,~k_2=k_3=l_1=0,~ l_2=-1,~ l_3=2$ and metric $g$ is isometric to Taub-NUT metric;

2. $k_1=1,~k_2=0,~k_3=-2,~ l_1=0,~ l_2=-1,~ l_3=0$ and metric $g$ is isometric to Eguchi-Hanson metric;

3. $k_1=1,~k_2=k_3=l_1=0,~ l_2=-1,~ l_3=-2$.}

Of course, if the set $(k_1,\dots,l_3)$ will give some metric $g_0$ then the set $(-k_i,-l_j)$ will give the same metric $g_0$, so we can exclude such cases.

\medskip
{\bf Theorem 2.}
{\it If metric $g$ is Einstein but not Ricci-flat and functions $A_1,~A_2$ of the metric satisfies $(1)$ and do not coincide then metric $g$ is isometric to Fubini-Study (or hyperbolic Fubini-Study) metric on $\mathbb{C}P^2$ (or on $SU(1,2)/S(U(1,1)\times U(1))$ respectively). For this case  $k_1=2,~k_2=0,~k_3=-1,~l_1=0,~l_2=1,~l_3=0$.}

We denote the ratio $\frac{A_1}{A_2}$ as a new variable $x$. If we assume that $x$ is a constant then $Ric_{00}=-6\frac{A_1''}{A_1}$, and $A_1(t)=\alpha t + \beta$ if one wants for $g$ to be Ricci-flat.
From the vanishing of components $Ric_{11}$ and $Ric_{22}$ one can deduce that $\alpha=1$, this corresponds to the flat Euclidean metric on $\mathbb{R}^4$. Next we will assume that $x$ is not a constant.

To prove first theorem it is sufficient to substitute expressions from $(1)$ into the components of the Ricci-tensor.
$$
Ric_{00}=- \frac {2}{{A_{{1}}}^{2}{x}^{6}}  \left( k_{{3}}{x}^{3}+k_{{2}}{x}^{2}-l_
{{3}}{x}^{2}+xk_{{1}}-l_{{2}}x-l_{{1}} \right) \cdot
$$
$$
\quad \quad \cdot \left(k_{{2}}{x}^{2}+ 2\,xk_{{1}}+2\,l_{
{2}}x+4\,l_{{1}} \right) .
$$

Then we have two cases.

a) $k_2=l_1=0$ and $k_1=-l_2$;

b) $k_3=l_1=0$, $k_2=l_3$ and $k_1=l_2$.

Substituting values from case a) to the $Ric_{11}$ one will get
$$
Ric_{11}=-4\,{\frac {{l_{{2}}}^{2}{x}^{3}+k_{{3}}l_{{3}}-{x}^{3}}{x{A_{{2}}}^{2}}}.
$$
So $l_2^2=1$ and $k_3l_3=0$. We can easily put $l_2=-1$. For $k_3=0$ one will have
$$
Ric_{22}=-2\,{\frac {{l_{{3}}}^{2}-4}{{A_{{2}}}^{2}}}
$$
and this are the cases 1 and 3 from the theorem 1.
For $l_3=0$ one will obtain case 2 from the theorem.

Case b) implies that
$$
Ric_{11}=-4\,{\frac {{l_{{2}}}^{2}{x}^{2}+2\,l_{{2}}l_{{3}}x+{l_{{3}}}^{2}-{x}^{2}}{{A_{{2}}}^{2}}},
$$
so one can put $l_3=0$ and $l_2=1$ into the $Ric_{22}$ and calculate that
$$
Ric_{22}=-8\,{\frac {{x}^{2}-1}{{A_{{2}}}^{2}}}.
$$
Thus case b) has no nontrivial solutions.

It is left to check that metrics of the cases 1 and 2 of the theorem are isometric to the previously known.
For the case 1 system $(1)$ turns to
$$
\left\{\begin{array}{ll}
A_1'=(\frac{A_1}{A_2})^2, \\
A_2'=-\frac{A_1}{A_2} +2.
\end{array}\right.
$$
For our choice of the basis $\{e^i\}$ Taub-NUT metric (see \cite{Taub} and \cite{NUT}) will take form
$$
\frac{r+m}{16(r-m)}dr^2 + \frac{m^2(r-m)}{r+m}e_1^2 + \frac{r^2-m^2}{4}(e_2^2+e_3^2). \eqno{(2)}
$$
Then
$$
A_1'=\frac{d A_1(t)}{dt}=\frac{d A_1(r)}{dr}\cdot \frac{dr}{dt} = \frac{m^2}{(r-m)^{\frac{1}{2}}(r+m)^{\frac{3}{2}}}\cdot \frac{4(r-m)^{\frac{1}{2}}}{(r+m)^{\frac{1}{2}}}=
$$
$$
\quad\quad\quad =\frac{4m^2}{(r+m)^2}=\frac{A_1^2}{A_2^2}.
$$
So the first equation is hold. Then
$$
A_2'=\frac{d A_2(r)}{dr}\cdot \frac{dr}{dt} = \frac{r}{2(r^2-m^2)^{\frac{1}{2}}}\cdot \frac{4(r-m)^{\frac{1}{2}}}{(r+m)^{\frac{1}{2}}}=\frac{2r}{r+m},
$$
$$
-\frac{A_1}{A_2} +2=-\frac{m(r-m)^{\frac{1}{2}}}{(r+m)^{\frac{1}{2}}} \cdot \frac{2}{(r^2-m^2)^{\frac{1}{2}}}+2=\frac{-2m}{r+m}+2=\frac{2r}{r+m}.
$$
And the second equation is true also.
One can easily check that Eguchi-Hanson metric \cite{EgHan}
$$
ds^2=[1-(a/r)^4]^{-1}dr^2 +r^2((e^2)^2+(e^3)^2)+r^2[1-(a/r)^4](e^1)^2
$$
is satisfied the system
$$
\left\{\begin{array}{ll}
A_1'=-\frac{A_1^2}{A_2^2}+2, \\
A_2'=\frac{A_1}{A_2}.
\end{array}\right.
$$
Which is exactly the case 2 up to a multiplication by $-1$.

To make the prove complete we must show that there is no other solutions of the system $(1)$ for $k_1=1,~k_2=k_3=l_1=0,~ l_2=-1,~ l_3=2$ except for the Taub-NUT. But this is simple conclusion of the counting free parameters of the system and metric. From the general theory of ODE it follows that arbitrary solution of $(1)$ depends on two parameters, namely $A_1(t_0)$ and $A_2(t_0)$. And for the Taub-NUT metric there are also two parameters: $m$ from $(2)$ and the time-shift $t_0$, i.e. one can always make a shift $t\rightarrow t+t_0$ without changing a metric $g$. The same is true for the Eguchi-Hanson metric. \textbf{Theorem 1} is proved.

\vskip0.5cm

To prove \textbf{Theorem 2} one should consider the expressions $\frac{Ric_{ii}}{g_{ii}}-\frac{Ric_{jj}}{g_{jj}}$ for $i\neq j \in\{0,1,2\}$. We emphasize that the components of the Ricci tensor were calculated for orthonormal frame $\{\varepsilon_0,\ldots, \varepsilon_3\}$, so actually one have to find solutions of following two equations
$$
Ric_{00}-Ric_{11}=0,\quad Ric_{00}-Ric_{22}=0,
$$
and then check if some of the components is actually a constant.

One can check using some mathematical programs or even by hands that
$$
Ric_{00}-Ric_{11}=\frac{4}{A_1A_2^6}\cdot (2l_1^2A_1^5+k_3l_3A_2^5 + Q(A_1,A_2)),
$$
where $Q(A_1,A_2)$ is a homogeneous polynomial of degree $5$ that does not contain monomials of the form $cA_1^5$ or $cA_2^5$. We suppose that functions $A_1$ and $A_2$ are linearly independent, so $l_1=0$. After one put $l_1=0$ he get
$$
Ric_{00}-Ric_{11}=\frac{4}{A_1A_2^4}\cdot ((l_2^2-1)A_1^3+ (k_1l_3+l_2l_3)A_1^2A_2 + k_2l_3A_1A_2^2 + k_3l_3A_2^3).
$$
Consider two cases: \emph{a)} $k_3=0$ and \emph{b)} $l_3=0$. Firstly we put $k_3=0$ and $l_2=-1$, then
$$
Ric_{00}-Ric_{11}=\frac{4l_3}{A_2^3}\cdot ((k_1-1)A_1+k_2A_2).
$$
If one will put $l_3=0$ then
$$
Ric_{00}-Ric_{22}=-\frac{2}{A_2^4}\cdot ((k_2^2+4)A_2^2+\ldots),
$$
and there are no solutions. Thus we put $k_1=1$ and $k_2=0$ and
$$
Ric_{00}-Ric_{22}=2(l_3^2-4)A_2^{-2}.
$$
This sets of parameters define cases 1 and 3 of the theorem 1, so we'll skip them.
If we assume that $l_3=0$ and $l_2=-1$ then $Ric_{00}-Ric_{11}=0$ and
$$
Ric_{00}-Ric_{22}=-\frac{2}{A_1A_2^4}\cdot((2k_1^2+2k_1-4)A_1^3 +k_2k_3A_2^3 +
$$
$$
\quad\quad\quad\quad +(2k_1k_3+k_2^2+4)A_1A_2^2 +(3k_1k_2+2k_2)A_1^2A_2).
$$
Then consider two cases a) $k_2=0$ and b) $k_3=0$. For a)
$$
Ric_{00}-Ric_{22}=-\frac{4}{A_2^4}\cdot((k_1^2+k_1-2)A_1^2+(k_1k_3+2)A_2^2),
$$
and one must put $k_1=1,~k_3=-2$ and this is a case 2 of the theorem 1 or $k_1=-2, k_3=1$. For b)
$$
Ric_{00}-Ric_{22}=-\frac{2}{A_2^4}\cdot((k_2^2+4)A_2^2+\ldots)
$$
and there is no solutions.
Now we want to multiply the set $\{k_1,\dots,l_3\}$ by $-1$ and to integrate the system
$$
\left\{\begin{array}{ll}
A_1'=2(\frac{A_1}{A_2})^2-1, \\
A_2'=\frac{A_1}{A_2}.
\end{array}\right.
$$
It is easy to verify that functions $A_1(t)=\frac{1}{2\alpha}\sin(2\alpha t)$ and $A_2(t)=\frac{1}{\alpha}\sin(\alpha t)$ are satisfy this system. For this solutions we have
$$
Ric_{00}=-2\frac{A_1''}{A_1}-4\frac{A_2''}{A_2}=-2[\frac{-2\alpha\sin(2\alpha t)}{\frac{1}{2\alpha}  \sin(2\alpha t)}  +2\frac{-\alpha\sin(\alpha t)}{\frac{1}{\alpha} \sin(\alpha t)}]=12\alpha^2,
$$
thus the found metric is actually Einstein.
For $\alpha=1$ one gets classical Fubini-Study metric, and for $\alpha=\mathfrak{i}$ one gets hyperbolic Fubini-Study with sectional curvature belongs to $[-4,-1]$, it is also well-known (\cite{KoNo}). The uniqueness of the Fubini-Study metric follows from the same arguments as in previous proof. \textbf{Theorem 2} is proved.

\section{Analysis of the case 3.}

We are not aware of any metric which functions would satisfy the system
$(1)$ for $k_1=1,~k_2=k_3=l_1=0,~ l_2=-1,~ l_3=-2$. The main problem of considered
system is that corresponding metric can not be defined on a complete
manifold: there is always a time moment $t_0$ such that
$A_1(t)\rightarrow + \infty$ and $A_2(t)\rightarrow 0$ as
$t\rightarrow t_0$.

Now we will find the explicit form of the metric. System $(1)$ turns to
$$
\left\{\begin{array}{ll}
A_1'=-(\frac{A_1}{A_2})^2, \\
A_2'=\frac{A_1}{A_2}+2.
\end{array}\right. \eqno{(3)}
$$
Denote the $A_1^{-1}$ as a new variable $\rho$. Then using the first equation from $(3)$ we can write
$$
d(A_1^{-1})=-A_1^{-2}\cdot \frac{\partial A_1}{\partial t}dt=A_2^{-2}dt=d\rho.
$$
The second equation from $(3)$ turns to $\dot{A_2}=\frac{A_2}{\rho}+2A_2^2$,
where $\dot{A_2}=\frac{\partial A_2}{\partial \rho}$. This latter equation can easily be integrated: $A_2(\rho)=\frac{\rho}{c^2-\rho^2}$, where $c^2$ is an integration constant. So metric $g$ takes the form
$$
\frac{\rho^4}{(c^2-\rho^2)^4}d\rho^2 +\frac{1}{\rho^2}(e^1)^2 +\frac{\rho^2}{(c^2-\rho^2)^2}((e^2)^2+(e^3)^2).\eqno{(4)}
$$
In this case dependence between $\rho$ and $t$ can be explicitly integrated:
$$
t= \int \frac{\rho^2}{(c^2-\rho^2)^2}d\rho=\frac{1}{4}\big( \frac{2\rho}{c^2-\rho^2}-\frac{1}{c}\ln \frac{\rho+c}{\rho-c}\big ),
$$
but it can not be explicitly inversed in elementary functions.

Next we want to understand the behavior of the functions $A_1(t),~A_2(t)$ as $t$ goes to infinity. Consider a function $B_2(t)=\alpha +\beta \ln t +2t$ and $B_1(t)=B_2(t)(B_2'(t)-2)=2\beta+\frac{\alpha\beta}{t}+\frac{\beta^2}{t}\ln{t}$ for some $\alpha,~\beta \in \mathbb{R}$. Then the second equation from $(3)$ is true identically and first one will imply that
$$
B_1'+\frac{B_1^2}{B_2^2}= \frac{\beta^2-\alpha\beta-\beta^2 \ln t}{t^2}+(\frac{\beta}{t})^2 \rightarrow 0.
$$
So the functions $B_1$ and $B_2$ are approximations of the $A_1$ and $A_2$ as $t\rightarrow \infty$.

We remind that a metric $dt^2 +\sum (A_i(t)e^i)^2$ is called Asymptotically Locally Conical (ALC) if it's functions $A_i(t)$ at infinity look like a linear functions, i.e. $\frac{A_i(t)}{a_it+b_i}\rightarrow 1$ as $t\rightarrow \infty$ for some constants $a_i,b_i$. And metric is called Asymptotically Locally Euclidean (ALE) if $a_i=1$. The right asymptotic at infinity can be useful for constructing compact examples using some gluing procedures (see \cite{Donaldson}).

For metric $(4)$ $A_1(t)\rightarrow const$ and $a_1=0$. We know that such behavior appears for metric with exceptional holonomy in dimension 7 and 8 (see \cite{BazMalk},\cite{BazBogo}). In dimension 4 Taub-NUT metric has the same asymptotic.

Now we will show that $A_1(t)\rightarrow + \infty$ and $A_2(t)\rightarrow 0$ as $t\rightarrow t_0$  for some time $t_0$. For this we consider the following functions
$$
C_2(t)=\gamma (t-t_0)^{\frac{1}{3}}  + \frac{9}{5}(t-t_0),
$$
$$
C_1(t)=C_2(t)(C_2'(t)-2)=  -\frac{(-5\gamma +3(t-t_0)^{\frac{2}{3}})(5\gamma(t-t_0)^{\frac{1}{3}}  +9(t-t_0))}{75(t-t_0)^{\frac{2}{3}}}.
$$
The second equation from $(3)$ is automatically satisfied as earlier. One can check that the difference between the left-hand side of the first equation from $(3)$ and the right-hand side is:
$$
C_1'+\frac{C_1^2}{C_2^2}=-\frac{8}{25}.
$$
But in this difference the first summand goes to $-\infty$ and the second one goes to $+\infty$, so their difference $-\frac{8}{25}$ is negligible. And we consider functions $C_1$ and $C_2$ as good approximations of the solutions $A_1$ and $A_2$ at time $t_0$. More accurate analysis can be done using the explicit dependence between $\rho$ and $t$ but we will omit it here.

\medskip
{\bf Remark}. Earlier in \cite{Malk} we presented flows that give standard metrics of constant curvature and the Eguchi-Hanson metric. We showed that for Eguchi-Hanson metric the right-hand side of the flow could be entirely expressed using only the Ricci tensor and some invariant functions of it like determinant. Also we showed that this flow can be established as a flow of the renormalized contact 1-forms:
$$
(\ast \psi)'=d\psi, \quad \mbox{for } \psi\in \{\varepsilon_1,\varepsilon_2,\varepsilon_3\} \eqno{(5)}
$$
It means that the system $(1)$ for the case 2 is equivalent to the
flow $(5)$. We tried to find some 'weighted flow' of the contact
1-forms that depends somehow on the parameters $(k_1,\dots,l_3)$.
But it turns out that, for example, system $(1)$ for Fubini-Study is
not equivalent to any flow of 1-forms that depends on natural
operations like Hodge star or exterior differential. So we pose the
question \emph{Is there a structure on $S^3$ (or $\mathbb{R}P^3$)
and evolution equation of it that is equivalent to the system
$(1)$?} Here by a 'structure' we mean some set of tensorial fields that can
generate a metric. Probably one should consider only the case when
$k_2=l_1=0$ because any of the described metrics satisfies this
condition.

\end{document}